\documentclass[12pt,twoside, a4paper]{article}
\usepackage{latexsym}
\usepackage[english]{babel}
\usepackage{t1enc}
\usepackage{mathrsfs}
\usepackage{ifthen}
\usepackage{url}
\usepackage{epsfig}
\usepackage{amsbsy}
\usepackage{extarrows}
\usepackage{amsfonts}
\usepackage{amsmath}
\usepackage{amssymb}
\usepackage{amsthm}
\usepackage{amsxtra}
\usepackage{bbm}
\usepackage{authblk} 
\usepackage[margin=2.5cm]{geometry}
\usepackage{color}

\usepackage{verbatim}

\usepackage{hyperref} % HAS TO BE THE LAST PACKAGE TO BE LOADED! - provides links inside the document!!!

\numberwithin{equation}{section}
\numberwithin{figure}{section}

\newtheoremstyle{thm-style-oskari}
{7pt}      % Space above
{7pt}      % Space below
{\itshape} % Body font
{}         % Indent amount (empty = no indent, \parindent = para indent)
{\scshape} % Thm head font
{.}        % Punctuation after thm head
{.5em}     % Space after thm head: " " = normal interword space; 
{}         % Thm head spec (can be left empty, meaning `normal')

\theoremstyle{thm-style-oskari}
    \newtheorem{theorem}{Theorem}[section]
    \newtheorem{proposition}[theorem]{Proposition}
    \newtheorem{corollary}[theorem]{Corollary}
    \newtheorem{lemma}[theorem]{Lemma}
    \newtheorem{definition}[theorem]{Definition} % howto make rm-style text inside definitions
    
    \newtheorem{convention}{Convention}[section]
    \newtheorem{remark}{Remark}[section]

\newenvironment{Proof}[1][Proof]{\begin{proof}[\sc{#1}]}{\end{proof}}

\newcommand{\Theorem}[2] {
        \begin{theorem} \label{thr:#1}
                #2
        \end{theorem}
        }

\newcommand{\Lemma}[2] {
        \begin{lemma} \label{lmm:#1}
                #2
        \end{lemma}
        }

\newcommand{\NCorollary}[2] {
        \begin{corollary}[#1] \label{crl:#1}
                #2
        \end{corollary}
        }

\newcommand{\bels}[2] {
        \begin{equation} \label{#1} \begin{split} 
                #2 
        \end{split} \end{equation}
        }

\newcommand{\mrm}[1] {\mathrm{#1}}
\newcommand{\mcl}[1] {\mathcal{#1}}

\newcommand{\brm}[1] {\boldsymbol{\mathrm{#1}}}

\newcommand{\wti}[1] {\widetilde{#1}}
\newcommand{\wht}[1] {\widehat{#1}}

\newcommand{\tE} {\mathbbm{E}}
\newcommand{\tP}  {\mathbbm{P}}

\newcommand{\abs}[1]{\lvert #1 \rvert}
\newcommand{\absb}[1]{\big\lvert #1 \big\rvert}
\newcommand{\absB}[1]{\Bigl\lvert #1 \Bigr\rvert}
\newcommand{\absbb}[1]{\biggl\lvert #1 \biggr\rvert}

\newcommand{\norm}[1]{\lVert #1 \rVert}
\newcommand{\normb}[1]{\big\lVert #1 \big\rVert}

\newcommand{\R} {\mathbb{R}}
\newcommand{\C} {{\mathbb{C}}}

\newcommand{\sett}[1] { \{ {#1} \} }

\newcommand{\setb}[1] { \bigl\{ {#1} \bigl\} }

\newcommand{\eps}{\varepsilon}

\newcommand{\msp}[1] {\mspace{#1 mu}} % put white space: valid input is integer
\newcommand{\1} {\mspace{1 mu}}
\newcommand{\2} {\mspace{2 mu}}

\newcommand{\mat}[1]{\begin{bmatrix} #1 \end{bmatrix}}

\newcommand{\cI} {\mathrm{i}}
\newcommand{\nE} {\mathrm{e}}

\renewcommand{\Im}{\mathrm{Im}}
\renewcommand{\Re}{\mathrm{Re}}

\newcommand{\tsfrac}[2] {{\textstyle \frac{#1}{#2}}}

\begin{document}
\title{Local semicircle law with imprimitive variance matrix}
\author[1]{Oskari Ajanki\thanks{Partially supported by SFB-TR 12 Grant of the German Research Council. On leave from Institute of Mathematics, University of Munich, \url{oskari.ajanki@ist.ac.at}}}
\author[1]{L\'aszl\'o Erd{\H o}s\thanks{Partially supported by SFB-TR 12 Grant of the German Research Council. On leave from Institute of Mathematics, University of Munich, \url{lerdos@ist.ac.at}}}
\author[1]{Torben Kr\"uger\thanks{Partially supported by SFB-TR 12 Grant of the German Research Council, \url{torben.krueger@ist.ac.at}}}
\affil[1]{IST Austria, Am Campus 1, Klosterneuburg A-3400}

\renewcommand\Authands{ and }

\maketitle
\thispagestyle{empty}
\begin{abstract}
We extend the proof of the local semicircle law for generalized Wigner matrices given in \cite{Erdos13-LoclScGenRMT} to the case when the matrix of variances has an eigenvalue $-1$. In particular, this result provides a short proof of  the optimal local Marchenko-Pastur law at the hard edge (i.e. around zero) for  sample covariance matrices $ \brm{X}^\ast \brm{X} $, where the variances of the entries of $\brm{X}$ may vary.
\end{abstract}

\section{Model and results}

The local semicircle law on the local distribution of eigenvalues of large Wigner matrices has been the basic technical input in the recent works on the Wigner-Dyson-Gaudin-Mehta universality (see \cite{EYBull} and references therein). The analysis was extended to generalized Wigner matrices \cite{Erdos11-BulkUnGenWig, Erdos13-LoclScGenRMT} but always practically assuming that the matrix of variances is primitive \footnote{A non-negative $ d \times d $-matrix $ \brm{M} $ is said to be primitive  (cf. Definition 4.1 of \cite{Berman94-NonNegMatrices}) if there  exists an integer $ k $ such that every element of $ \brm{M}^k $ is strictly positive.}, in particular $ -1 $ is not in its spectrum. 
This assumption naturally holds for random band matrices that were the main motivation to generalize Wigner matrices in  \cite{Erdos11-BulkUnGenWig}. However, some important matrices with a certain block structure do not satisfy this condition. Most notable example is the $2N\times 2N$ matrix
\bels{HXX}{
\brm{H} = \mat{ \,\brm{0} & \brm{X^*} \\ \2\brm{X} & \!\brm{0}\, }
}
where the $ N \times N $ matrix $ \brm{X} $ has independent entries. The matrix $ \brm{H}$ is the linearization of the of the sample covariance matrix $\brm{X^*}\brm{X}$. In this paper we show how to remove the primitivity assumption in \cite{Erdos13-LoclScGenRMT}.

We consider generalized $ N\times N $ hermitian or symmetric Wigner matrix $ \brm{H} = (h_{ij})_{i,j=1}^N $ with independent entries (up to the symmetry constraint $\brm{H}=\brm{H}^*$) such that
\bels{H is mean zero and s_ij finite}{
\tE\, h_{ij} \,=\,0\,,
\qquad\text{and}\qquad
s_{ij} := \tE\,\abs{h_{ij}}^2 < \infty
\,.
}
We assume that all moments are bounded in the sense that,
\bels{finite moments}{
\tE\, \absbb{\frac{h_{\1i\1j}}{s_{ij}^{1/2}\!}}^p <\,  C_p\,,
\qquad
\forall\,p < \infty
\,,
\msp{-50}
}
with constants $ C_p $ independent of $ N $.
In order to avoid unnecessary clutter we have suppressed the $ N $-dependence in the notations, e.g., we use $ \brm{H} $ and $ \brm{S}$ to refer to the sequences of matrices $ \brm{H}^{(N)} = (h_{ij}^{(N)})_{i,j=1}^N $ and $ \brm{S}^{(N)} = (s^{(N)}_{ij})_{i,j=1}^N $, respectively.

Besides the natural constraints, $ \brm{S}^{\mrm{T}} = \brm{S} $, $ s_{ij} \ge 0 $,
 we make the following additional assumptions on the variance matrix: 
\begin{itemize}
\item[(A1)] 
{\bf Boundedness:} There exists a sequence $  N^\delta \leq M = M_N \leq N $, with $ \delta > 0 $, such that  
\bels{bounded}{
0 \,\leq\, s_{ij} \,\leq\, M^{-1}\,;
}
\item[(A2)]
{\bf Constant row sums:} $ \brm{S} $ is (double) stochastic:
\bels{constant}{
\sum_{j=1}^N s_{ij} \,=\, 1
\,,\qquad
\forall\,i=1,\dots,N\,;  
}
\item[(A3)]
{\bf Isolated extremum eigenvalues:} There exists ($ N $-independent) constant $ 0 < \rho < 1 $,  such that 
\bels{isolated extremum eigenvalues}{ 
\mrm{Spec}(\brm{S}) \subset \sett{-1}\cup [-\rho,\2\rho\2] \cup \sett{+ 1}\,. 
}
\end{itemize}
This setup is similar to that in \cite{Erdos13-LoclScGenRMT}, except that here we explicitly allow $ -1 $ in the spectrum of $\brm{S} $. This allows us to consider  $ \brm{S} $ which contain imprimitive irreducible components. The results in \cite{Erdos13-LoclScGenRMT} practically excluded this case since the estimates became unstable, see Section 7 of \cite{Erdos13-LoclScGenRMT}. The main observation of this paper is that this instability is not present.

The relaxation of the irreducibility condition is elementary algebra (cf. Lemma \ref{lmm:Decomposition of S} below), and this extension was already mentioned in \cite{Erdos13-LoclScGenRMT}. However, the inclusion  of $ -1 $'s in the spectrum of $\brm{S}$ requires a new  algebraic identity that is stated as Lemma \ref{lmm:symm2} below. We will show here how to incorporate this identity into the proof given in \cite{Erdos13-LoclScGenRMT} with minor modifications.

The condition (A2) guarantees that the diagonal elements of the  resolvent matrix,
\bels{}{
\msp{100}\brm{G}(z) \,:=\, 
\frac{1}{\brm{H}-z}
\,,\qquad
z := E + \cI\1\eta\,,\quad E \in \R\,,\;\eta > 0
\,,
}
converge towards the Stieltjes transform 
\bels{}{
m(z) \,=\, 
\frac{-\1z \,+\2 \sqrt{z^2-4\,}}{2}
}
of the Wigner semicircle law, $\varrho(x) = (2\1\pi)^{-1}\sqrt{\max\sett{4-x^2,0}}$, as $ N $ approaches infinity.

In order to state this main result, we recall the concept of stochastic domination (Definition 2.1 in \cite{Erdos13-LoclScGenRMT}).
We say that  a (sequence) of random variables $ X=X^{(N)} $ is stochastically dominated by another (sequence) of random variables $ Y=Y^{(N)} $, in notation $ X \prec Y$, if for any $\eps, D>0$ there exists $ N_0 = N_0(\eps, D)< \infty $ such that
\bels{}{
\tP\setb{X > N^\eps Y}\;\leq\; N^{-D}
\qquad \forall\,N\ge N_0 
\,.
}
If $X, Y$ depend on some other parameters (like $z$ or labels like $i, j$), then the definition is always taken uniform in these parameters (i.e. $N_0$ depends only on $\eps, D>0$). The notation $ X = \mcl{O}_\prec(Y)  $ means same as $ \abs{X} \prec Y $. 

\Theorem{Main Result}{
Suppose $ \brm{S} $ satisfies the assumptions  (A1)--(A3), and denote
\bels{def of spectral domain D(gamma)}{
\mathbb{D}(\gamma) 
\,=\, 
\mathbb{D}^{(N)}(\gamma) 
\;:=\, \setb{ z  : \; \abs{z}\leq 10\,,\; \Im\,z\ge M_N^{-1+\gamma}}\,,\qquad\gamma > 0\,.
}
Then for any fixed $\gamma>0$, the local estimates
\begin{subequations}
\label{main results}
\begin{align}
\label{x-lsc}
\max_{i,\1j=1}^N\absb{G_{ij}(z) -m(z)\1\delta_{ij}}
\;\,&\prec\;
\sqrt{\frac{\Im \,m(z)}{M\eta}} +\frac{1}{M\eta}
\\
\label{lsc}
\absB{\frac{1}{N}\mrm{Tr}\,  \brm{G}(z) - m(z)}  \;&\prec\;\, \frac{1}{M\eta} 
\,,
%}
\end{align}
\end{subequations}
apply uniformly in $ z = E + \cI\1\eta \,\in\, \mathbb{D}(\gamma) $.

Moreover, outside the spectrum \eqref{lsc} can be strengthened by introducing the distance, $ \kappa := \max\sett{\2\abs{E-2},\abs{E+2}\2} $, of $ E $ from the spectral edges:
\bels{lscout}{
\absB{\frac{1}{N}\mrm{Tr}\,\brm{G}(z) - m(z)}  
\;\prec\; 
\frac{1}{M(\kappa+\eta) } \,+\,
\frac{1}{(M\eta)^2\sqrt{\kappa+\eta}}
\;.
}
This estimate is also uniform in $ z =  E+\cI\1\eta \in \mathbb{D}(\gamma) $ as long as the constraints $ \abs{E}\ge 2 $ and $ \eta\sqrt{\kappa+\eta} \ge M^{-1+\gamma} $ are satisfied.
}

This theorem is a generalization of the following result.

\begin{theorem}[\cite{Erdos13-LoclScGenRMT}] 
\label{thr:OldResult}
Assume that $ \brm{S} $ satisfies (A1)--(A2), but (A3) is strengthened to 
$ \mrm{Spec}(\brm{S}) \subset [-\rho,\rho\1] \cup \sett{+1} $.
Then conclusions of Theorem \ref{thr:Main Result} hold true.
\end{theorem}

Theorem \ref{thr:OldResult} is a special case of the more general Theorem 2.3 of \cite{Erdos13-LoclScGenRMT} which also covers the cases $ \mrm{Spec}(\brm{S}) \subset [-\rho^{(N)}_-\, ,\, \rho_+^{(N)} \1] \cup \sett{+1} $ where the spectral gaps $1-\rho_\pm^{(N)} $ may close at certain rates as $ N \to \infty $. 
However, the estimates near $ E = 0 $ in \cite{Erdos13-LoclScGenRMT} deteriorated if the smallest eigenvalue approached  $ -1 $, and in particular $ -1 $ was not allowed belong to the spectrum. The condition (A3) rules out closing of the upper gap, and hence the spectral domain $ \wti{\mathbb{S}}(\gamma) $, defined by formulas (2.14) and (2.17) in \cite{Erdos13-LoclScGenRMT}, has been replaced here by the simpler set $ \mathbb{D}(\gamma) $. 
It is straighforward to extend Theorem~\ref{thr:Main Result} to the entire set  $ \wti{\mathbb{S}}(\gamma) $ in the spirit of Theorem 2.3 in \cite{Erdos13-LoclScGenRMT} but for brevity of this note we refrain from doing so. 

Theorem~\ref{thr:Main Result} directly implies a rigidity result for the increasingly ordered eigenvalues $ (\lambda_\alpha)_{\alpha=1}^N $  of $\brm{H} $ in terms of the  $N$-th quantiles $ (\gamma_\alpha)_{\alpha=1}^N$ of the semicircle density:
\bels{rig}{
\abs{ \2\lambda_\alpha -\2\gamma_\alpha\1}
\;\prec\;
\frac{1}{M}
\Big(\frac{N}{\widehat\alpha}\Big)^{1/3}
\,,
\quad\text{when}\quad
\wht{\alpha} := \min\sett{\alpha, N+1-\alpha} \,\ge\, N M^{-1+\eps}
\,,
}
with $ \eps > 0 $ arbitrary. 
See Theorem 7.6 in \cite{Erdos13-LoclScGenRMT} for a proof in a more general setup and for the estimates on the extreme eigenvalues.

We remark that there have been  many  results on local semicircle laws prior to  \cite{Erdos13-LoclScGenRMT}, in fact most methods used in \cite{Erdos13-LoclScGenRMT} stem from \cite{Erdos11-BulkUnGenWig, Erdos11-UnivGenWigBernoulli, Erdos12-RigGenWig}. See  \cite{Erdos13-LoclScGenRMT} for a complete account of the history and for the most concise general proof.

Finally, we mention a simple application. The eigenvalues of $ \brm{H} $ in \eqref{HXX} generically come in pairs, $\pm\lambda $, (see \eqref{simple eigenpairs of H} below) and their squares $\lambda^2 $ are the  eigenvalues of the sample covariance matrices $\brm{X}\1\brm{X}^*$ and $\brm{X}^*\brm{X}$.
We assume that the elements of the square matrix $ \brm{X} $ are independent, centred and their variances are chosen such that $\brm{S}$ and $\brm{H} $ satisfy \eqref{finite moments}--\eqref{isolated extremum eigenvalues} (note that $-1$ is an eigenvalue of $\brm{S}$). 
Under these conditions, Theorem~\ref{thr:Main Result} can be directly used to estimate the  resolvent matrix elements and the trace of $\brm{X}^*\brm{X}$. Indeed, by applying the Schur formula to the $ N \times N $-block decomposition
\[
\brm{G}(z) 
\,=\, 
\frac{1}{\brm{H}-z}
\,=\,
\mat{
-z & \brm{X}^\ast\, 
\\
\brm{X}\! & -z\,}^{-1}
\!=:\, 
\mat{ 
\1\brm{G}_{11} &  \brm{G}_{12} 
\\  
\1\brm{G}_{21} & \brm{G}_{22}
}
\,,
\]
we see that the blocks on the diagonal equal:
\[
\brm{G}_{11}(z) \,=\, 
\frac{z}{\brm{X^*}\brm{X} -z^2}, 
\qquad  
\brm{G}_{22}(z) \,=\, 
\frac{z}{\brm{X}\brm{X^*} -z^2}
\,.
\]
 Thus Theorem~\ref{thr:Main Result} implies that the local Marchenko-Pastur law holds in the critical ``hard-edge'' case, when the limiting density\\
\bels{MP-density}{
\varrho_{\mrm{MP}}(x) \,:=\, \frac{1}{2\pi}\sqrt{\frac{\max\sett{4-x,0}}{x}}
\,,
}
is singular at the origin. In fact $ \varrho_{\mrm{MP}}(x) = x^{-1/2}\varrho(x^{1/2})$, $x>0$,  where $ \varrho $ is the Wigner semicircle density.
By denoting the Stieltjes transform of the Marchenko-Pastur law \eqref{MP-density} by $ m_{\mrm{MP}} $ and and writing $ w := z^2 $, an elementary calculation from Theorem~\ref{thr:Main Result} yields the following result.

\NCorollary{Local Marchenko-Pastur law at the hard edge}{
Under the conditions on $\brm{X}$ above, we have for any fixed $ \gamma > 0 $,
\begin{subequations}
\label{}
\begin{align} 
\max_{i,j=1}^N\,\absbb{\Big( \frac{1}{\brm{X}^*\brm{X}- w\1} \Big)_{ij} -\, m_{\mrm{MP}}(w)\1 \delta_{ij}\,} 
\;&\prec\: 
\sqrt{\frac{\Im \,m_{\mrm{MP}}(w)}{M\2\Im\,w}} \,+\frac{1}{M\2\Im\,w}
\\
\absbb{\,\frac{1}{N}\,
\mrm{Tr} \, \frac{1}{\brm{X}^*\brm{X}- w\1} \,-\, m_{\mrm{MP}}(w)\,} 
\;&\prec\;\, 
\frac{1}{M\2\Im\,w}
\,,
\end{align}
\end{subequations}
uniformly in $ w \in \C $ satisfying $ \abs{w} \leq 100 $ and $\Im\,w \ge 
\sqrt{\abs{\Re\,w}\,}\,M^{-1+\gamma}$. 
}
The estimate outside of the spectrum \eqref{lscout} and the rigidity bound \eqref{rig} can also be directly translated to the similar statements for the sample covariance matrices.

We remark that local Marchenko-Pastur law on the smallest local scale was first proven in \cite{ESYY} away from the critical case. The hard-edge case was independently  considered in \cite{CMS} and in  \cite{Bourgade13-LocCircRMT}, the latter providing an optimal error bound. Both works dealt with the case when the variances  $ \tE\, \abs{x_{ij}}^2$ are constant, the above corollary extends the result to the case of non-constant variances.

\medskip

{\it Acknowledgement.} The authors are grateful to Ofer Zeitouni who  pointed out  the importance of removing the primitivity condition on $\brm{S}$ for various applications. 
An alternative albeit somewhat weaker extension of  \cite{Erdos13-LoclScGenRMT} to treat the case $ -1 \in \mrm{Spec}(\brm{S}) $ was given in \cite{VZ}.

\section{Two algebraic lemmas}

Let us define for arbitrary square matrices $ \brm{M}_1,\dots,\brm{M}_k $, the diagonal block matrix by:
\bels{def of block matrix brmD}{
\brm{D}(\brm{M}_1,\dots,\brm{M}_k)
\;:=\,
\mat{
\,\brm{M}_1\! & \brm{0} & \cdots & \!\brm{0} 
\\
\brm{0} & \brm{M}_2\! & \cdots & \!\brm{0}
\\
\vdots & \vdots & \ddots & \!\vdots 
\\
\brm{0} & \brm{0} & \cdots & \brm{M}_k
}
\,.
}

General algebraic results for non-negative matrices yield the following decomposition when applied to $ \brm{S} $ satisfying the assumptions of Theorem \ref{thr:Main Result}. 

\Lemma{Decomposition of S}{
Suppose that $ \brm{S} $ is a symmetric (double) stochastic matrix that satisfies 
conditions (A1)--(A3). 
Then, after an appropriate permutation $ \brm{P} $ of the indices, $ \brm{S} $ has a  block structure
\bels{irreducible decomposition of S}{
\brm{S} 
\,=\,
\brm{P}
\2
\brm{D}\bigl(\2\brm{S}_1,\dots,\brm{S}_p\1,\2\wti{\brm{S}}_1,\dots,\wti{\brm{S}}_q\bigr)
\,
\brm{P}^{-1}
\,,
}
where $ \brm{S}_\alpha $, $ 1 \leq \alpha \leq p $, and $ \wti{\brm{S}}_\beta $, $ 1 \leq \beta \leq q $, are irreducible doubly stochastic matrices with some $ p,q $.
The spectrums of the blocks
\[ 
\wti{\brm{S}}_\beta = (\tilde{s}_{\beta;\1ij})_{i,j=1}^{\tilde{d}_\beta}
\,,
\quad 1 \leq \beta \leq q\,,
\]
satisfy $ \mrm{Spec}(\wti{\brm{S}}_\beta)\subset
 [-\rho,\2\rho\1]\cup\sett{+1} $.
The blocks $ \brm{S}_\alpha $ have both $ +1 $ and $-1$ as simple eigenvalues, and they have the  structure
\bels{irreducible block of S with period 2}{
\qquad 
\brm{S}_\alpha =  \mat{\brm{0} & \brm{A}_\alpha^{\!\mrm{T}} \\ \brm{A}_\alpha\! & \brm{0}\;}
\,,\qquad
\brm{A}_\alpha = (a_{\alpha;\1ij})_{i,j=1}^{d_\alpha}\,,
\quad 1 \leq \alpha \leq p\,,\msp{-30}
}
where both the rows and the columns of the (generally) non-symmetric matrices $\brm{A}_\beta $  sum to one, i.e., $ \sum_i a_{\beta;\1ij}= \sum_j a_{\beta;\1ij}=1 $.
The matrix elements are bounded,
\[
\tilde{s}_{\beta;\1ij}\,,\,a_{\beta;\1ij}  \,\leq\, \frac{1}{M}
\,,
\]
and the dimensions satisfy  $ d_\alpha,\,\tilde{d}_\beta \ge M $, and
\bels{}{
2\1d_1 \,+\, \dots \,+\,2\1d_p \,+\, \tilde{d}_1 \,+\, \dots \,+\,\tilde{d}_q \;=\; N
\,. 
}
}
\begin{Proof}
Irreducible components of $ \brm{S} $ can be permuted into diagonal square blocks \eqref{irreducible decomposition of S} and the properties (A1)--(A3) are preserved under relabelling. 
In particular, the constant row and column sums for $ \brm{S} $, as well as the bound $ 0 \leq s_{ij} \leq M^{-1} $, translate directly to analogous bounds for $ \wti{S}_\beta $'s and $ \brm{A}_\alpha $'s, and this in turn implies $ d_\alpha, \tilde{d}_\beta \ge M $.

The structure of the block decomposition of
$ \brm{S}_\alpha $  \eqref{irreducible block of S with period 2} follows from the general theory of non-negative irreducible matrices $ \brm{M} = (m_{ij})_{i,j=1}^d $, $ m_{ij} \ge 0 $, e.g. Theorem 2.20 of \cite{Berman94-NonNegMatrices}:
If $ \brm{M}$ has $ k $-eigenvalues on its spectral circle, $ \sett{z \in \C: \abs{z} = r }$, then those eigenvalues are precisely the $ k $ complex roots of $ r^2 $, i.e., they equal $ r\1\nE^{\cI\12\pi j/k}$,  $ 1 \leq j\leq k $. Moreover, the matrix $\brm{M}$ has the block representation: 
\bels{general decomposition of non-negative matrix M}{
\brm{M} 
\;=\; 
\mat{
\brm{0} & \brm{D}\bigl(\1\brm{M}_1,\dots,\brm{M}_{k-1}\bigr)
\\
\,\brm{M}_k\msp{-7} & \brm{0}
} 
}
for some matrices $ \brm{M}_1, \ldots, \brm{M}_k $.
The dimensions of $ \brm{M}_j $'s are such that the zero blocks along the diagonal (not visible in \eqref{general decomposition of non-negative matrix M}) are square, so that the rows of $ \brm{M}_j $ and the columns of $ \brm{M}_{j+1} $ have the same dimensions for each $j=1,2,\ldots,k $ if one identifies $ \brm{M}_{k+1} := \brm{M}_1 $.

Applying the decomposition \eqref{general decomposition of non-negative matrix M} to $ \brm{M} := \brm{S}_\alpha $ yields the representation \eqref{irreducible block of S with period 2} since the symmetry $ \brm{S}_\alpha^{\mrm{T}} = \brm{S}_\alpha $ implies $ 1 \leq k \leq 2 $, while $ -1 \in \mrm{Spec}(\brm{S}_\alpha) $ excludes the case $ k \neq 1 $.
\end{Proof}

The random matrix $ \brm{H} $ inherits the structure \eqref{irreducible decomposition of S} of $ \brm{S} $ through \eqref{H is mean zero and s_ij finite}:  
\[
\brm{H}
\,=\,
\brm{P}
\2
\brm{D}\bigl(\2\brm{H}_1,\dots,\brm{H}_p\1,\2\wti{\brm{H}}_1,\dots,\wti{\brm{H}}_q\bigr)
\,
\brm{P}^{-1}
\,.
\]
Here  $\brm{H}_\alpha $ and $ \wti{\brm{H}}_\beta $ are independent generalised Wigner matrices satisfying 
$ \tE\,\abs{h_{\alpha;ij}}^2 = s_{\alpha;ij} $ and  $ \tE\,\abs{\tilde{h}_{\beta;ij}}^2 = \tilde{s}_{\beta;ij} $, respectively. 
This decomposition means that it suffices to prove Theorem \ref{thr:Main Result} for the irreducible components separately.
The components $ \wti{\brm{H}}_\beta $ are already covered by Theorem \ref{thr:OldResult}.
Hence, dropping the indices $ \alpha \ge 1 $, we are left to prove Theorem \ref{thr:Main Result} in the case
\bels{blocks of H and S}{
\brm{H} \,=\,  \mat{\;\brm{0} & \,\brm{X}^* \\ \brm{X}\! & \brm{0}\,}
\qquad\text{and}\qquad
\brm{S} \,=\,  \mat{\;\brm{0} & \,\brm{A}^{\!\mrm{T}} \\ \,\brm{A}\! & \brm{0}\,}
\,,
}
where $ \brm{S} $ is irreducible, and the entries $ x_{ij} $ of the square matrices $\brm{X}$ are independent, and satisfy $ \tE\,\abs{x_{ij}}^2 = a_{ij} $. 
For the sake of convenience, we also redefine $ N $ to be equal to the 
dimension of $ \brm{A} $, so that $ \brm{H} $ and $ \brm{S} $ are
 $ 2N \times 2N$ matrices.

Using the special structure \eqref{blocks of H and S} it follows that if $ \lambda \in \mrm{Spec}(\brm{H}) $ then also $ -\lambda \in \mrm{Spec}(\brm{H}) $, and the corresponding eigenvectors are related in the following simple way:
\bels{simple eigenpairs of H}{
\brm{H}\mat{\brm{u}\\ \pm\brm{w}} = \pm\,\lambda \mat{\brm{u}\\ \pm\brm{w}} 
\,.
}
In particular, the same reasoning can be applied to the $ \pm 1 $ eigenvalues of $ \brm{S } $. Moreover, since $ \brm{S} $ is double stochastic and irreducible, the  eigenvectors of $ \brm{S} $ belonging to the non-degenerate eigenvalues $ \pm 1$,   equal
\bels{def of e and f}{
\brm{e} \,&:=\, \tsfrac{1}{\sqrt{2N}}(\11,\dots,1,1,\dots,1)
\\
\brm{f} \,&:=\, \tsfrac{1}{\sqrt{2N}}(\11,\ldots, 1, -1, -1, \ldots, -1)
\,,
}
so that $ \brm{S}\1\brm{e} = \brm{e} $ and $ \brm{S}\1\brm{f} = -\1\brm{f}$.

Let us denote the complex inner product between vectors $ \brm{a},\brm{b} \in \C^{2N} $ by $ (\brm{a},\brm{b}) = \sum_i \overline{a_i}\,b_i $.     
Combining the symmetries  \eqref{simple eigenpairs of H} and \eqref{def of e and f} of $ \brm{H} $ and $ \brm{S} $ yields 
the following  algebraic identity.

\Lemma{symm2}{
Let $ \brm{H} $ be a $ 2N $-dimensional self-adjoint matrix that has the block structure \eqref{blocks of H and S} but is otherwise arbitrary, e.g., \eqref{H is mean zero and s_ij finite} is not assumed. 
Then the Green function $ \brm{G}(z) := (\brm{H}-z)^{-1} $ is orthogonal to the $ -1 $ eigenspace of $ \brm{S} $,
\bels{id2}{
\bigl(\1\brm{f}\1,\1\mrm{diag}(\brm{G})\bigr) \,=\, 0
\,,
} 
where $ \mrm{diag}(\brm{G}) = (G_{11},\dots,G_{2N,2N}) $ and $ \brm{f} $ is defined in \eqref{def of e and f}.
}
\begin{Proof}
Let us additionally assume that $ 0 \notin \mrm{Spec}(\brm{H}) $, and, that besides the pairs \eqref{simple eigenpairs of H}, there are no further degeneracies. 
Suppose
\[
\brm{v} := \mat{\brm{u}\\ \brm{w}}
\qquad\text{and}\qquad
\wti{\brm{v}} := \mat{\brm{u}\!\\ -\brm{w}\,}
\,,\qquad\text{with}\qquad \brm{u},\brm{w} \in \C^N\,,\msp{-100}
\]
are the eigenvectors corresponding to the eigenvalues $ \lambda $ and $ - \lambda $ in \eqref{simple eigenpairs of H}, respectively.
Since $ \lambda \neq 0 $ (so that $ -\lambda \neq \lambda $) these eigenvectors are orthogonal
\[
0\,=\, (\brm{v},\wti{\brm{v}}) \,=\, 
 \sum_{i=1}^N\, \abs{u_i}^2 - \sum_{k=1}^N \abs{w_k}^2
\,,
\]
i.e., the first and the second blocks are balanced: $ \norm{\brm{u}}_2 = \norm{\brm{w}}_2 $. 

Now, let $ \brm{v}^{(\alpha)} $, $ \alpha =1,\dots, 2N $, be the $ 2N $ eigenvectors of $ H $, and let $ \brm{u}^{(\alpha)} $ and $ \brm{w}^{(\alpha)} $ contain the first and the last $ N $-components of $ \brm{v}^{(\alpha)} $. 
Then combining the spectral theorem,
\[
G_{kk}(z) \,=\, \sum_{\alpha=1}^{2N} \frac{ |v_k^{(\alpha)}|^2 }{\lambda_\alpha-z}\,,
\]
with the balancing conditions, $ \norm{\brm{u}^{(\alpha)}}_2 = \norm{\brm{w}^{(\alpha)}}_2 $,  yields 
\bels{balancing conditioning for G}{
\sum_{k=1}^N G_{kk} 
\,=\, 
\sum_{\alpha=1}^{2N} \frac{\norm{\brm{u}^{(\alpha)}}_2^2}{\lambda_\alpha-z} 
\,=\,
\sum_{\alpha=1}^{2N} \frac{\norm{\brm{w}^{(\alpha)}}_2^2}{\lambda_\alpha-z}
\,=\,
\sum_{k=N+1}^{2N} G_{kk}
\,,
}
which is equivalent to \eqref{id2}. 
Finally, by using a basic continuity argument, one sees that \eqref{balancing conditioning for G} must apply also without the extra assumptions concerning the degeneracies and the exclusion of $ 0 $ from the spectrum of $ \brm{H}$.
\end{Proof}

\section{Translating proof of Theorem \ref{thr:OldResult} to cover  case \eqref{blocks of H and S}}

With the identity \eqref{id2} at hand we may translate the proof of Theorem \ref{thr:OldResult} from \cite{Erdos13-LoclScGenRMT} to the setting \eqref{blocks of H and S} without significant changes. 
In order to see this, we recall that the $ -1 $ eigenvalue of $ \brm{S} $ enters the proofs in \cite{Erdos13-LoclScGenRMT} only when one needs to bound the inverse of the operator $ 1-m^2\brm{S} $. 
However, using the identity \eqref{id2} one can  show that in all such cases it suffices to restrict the analysis to the orthogonal complement of the eigendirection $ \brm{f} $ corresponding to the eigenvalue $-1$.
The lower spectral gap assumption \eqref{isolated extremum eigenvalues} above $ -1 $ then guarantees that  $ (1-m^2\brm{S})^{-1} $  stays uniformly bounded even when $ m(z)^2 $ becomes close to $ -1 $ (equivalently, $z\approx 0$), i.e.,
\[
(1-m(z)^2\brm{S})^{-1} \,\approx\; (1 + \brm{S}\1)^{-1}\,,
\qquad\text{for}\quad z \approx 0\,.
\] 
Since the 'bad direction' $ \brm{f} $ will not play any role in the analysis, the only necessary modification of \cite{Erdos13-LoclScGenRMT} in the end,
 is to replace the operator norm $ \wti{\Gamma}(z) $ (cf. equation (2.11) in \cite{Erdos13-LoclScGenRMT}) of $ (1-m^2\brm{S})^{-1} $ in $ \brm{e}^\perp $, by the analogous norm in the orthogonal complement of both $ \brm{e} $ and $ \brm{f} $: 
\bels{}{
\wht{\Gamma}(z) \,:=\, \normb{ (1-m(z)^2\brm{S}\1)^{-1}|_{ \mrm{span}\sett{\brm{e}\1,\2 \brm{f}}^\perp}}_{\ell^\infty \to \ell^\infty}
\,.
}
The estimate (A.3) of \cite{Erdos13-LoclScGenRMT},
\[
\wht{\Gamma}(z) \,\leq\, C(\rho) \log N
\]
remains valid since the operator norm of $ (1-m^2\brm{S})^{-1} $ from $ \ell^{\12} $ to itself is bounded by $1/(1-\rho)$ in the complement of $  \mrm{span}\sett{\brm{e}\1,\2 \brm{f}} $. Here the logarithm comes from the fact that the $\ell^{\1\infty}$-norm is bigger by this factor over the $\ell^{\12}$-norm (cf. p. 46 of \cite{Erdos13-LoclScGenRMT}).

It remains to demonstrate why the inversion of $ 1-m^2\1\brm{S} $ can be always be restricted to the orthogonal complement of $ \brm{f}$. 
This inversion was used to bound the random fluctuations of the diagonal resolvent elements,
\bels{def of v}{
v_i \,:\,= G_{ii}-m_i 
}
in terms of the small random error terms $ \Upsilon_i = \mcl{O}_{\!\prec}(N^{-c}) $ appearing the self-consistent vector equation (cf. (5.9) in \cite{Erdos13-LoclScGenRMT}):
\bels{SCE for errors v}{
- \sum_k s_{ik} v_k \,+\, \Upsilon_i \;=\; \frac{1}{m+v_i} \,-\,\frac{1}{m}
\,.
}
Under the assumption, $ \abs{v_i} \prec \Lambda \prec N^{-c} $, with some control parameter $\Lambda$,  and 
using  $ \abs{m} \sim 1 $, \eqref{SCE for errors v}, takes the form 
\bels{linearised SCE for v}{
(1-m^2\brm{S}\1)\1\brm{v} \,=\, \mcl{O}_\prec(\norm{\brm{\Upsilon}}_\infty+\Lambda^2)  
\,.
}
Writing \eqref{def of v}  as 
\[
\brm{v} = \mrm{diag}(\brm{G}) - \sqrt{2N}\1m\2\brm{e} 
\,,
\]
recalling  $ (\brm{f},\brm{e}) = 0 $, and then applying Lemma \ref{lmm:symm2} yields:
\bels{f-projection of v = 0}{
(\1\brm{f},\brm{v}) \,=\, (\1\brm{f},\1\mrm{diag}(\brm{G})\1) \,=\,0 
\,.
}
The identity \eqref{f-projection of v = 0} shows that inversion of $ 1-m^2\brm{S} $ can be indeed restricted to the complement of $ \brm{f} $ in the case of \eqref{linearised SCE for v}.

The inverse of $ 1-m^2\brm{S} $ becomes unbounded also in the direction $ \brm{e} $ when $ m^2 \approx 1 $. 
However, unlike with direction $ \brm{f} $, the inversions of  $ 1-m^2\brm{S} $ can not be straightforwardly restricted to the complement of $ \brm{e} $, since the average of $ \brm{v} $, 
\bels{}{
[\brm{v}] \,:=\, \frac{1}{2N}\sum_{i=1}^{2N}v_i \,=\, \frac{(\brm{e},\brm{v})}{\!\sqrt{2N\,}} 
\,,
}
is not small.
For this reason the critical part $ [\brm{v}] $ was treated separately from the remainder, $ \brm{v} - (\brm{e},\brm{v})\2 \brm{e} \in \mrm{span}\sett{\brm{e},\brm{f}}^\perp $ in a more precise second order scalar equation in \cite{Erdos13-LoclScGenRMT}. 
The remainder part satisfies a linearised vector equation for which one needs to again invert $ 1-m^2\brm{S} $. We will now demonstrate that also in this case the component $ \brm{f} $ is not present due to Lemma \ref{lmm:symm2}. 
Indeed, in order to get from (6.19) to (6.20) in \cite{Erdos13-LoclScGenRMT} one applies the fluctuation averaging estimate (4.14) (with the choice $ t_{ij} = s_{ij}) $ to bound the remainder $ \brm{v} - (\brm{e},\brm{v})\2 \brm{e} $.
The crucial steps appear in the proof of (4.14) located at the end of the proof of Theorem 4.7 on p. 54 of \cite{Erdos13-LoclScGenRMT}, where a bound for
\bels{}{
w_a := \sum_i t_{ai}\1(v_i- [v])
\,, 
}
is derived from a linearised self-consistent equation
\bels{SCE for FA}{
\sum_i t_{ai}(v_i-[\brm{v}]) \,&=\; m^2 \sum_{b,j} s_{ab}\1t_{bj} (v_j-[\brm{v}]) \;+\; \mcl{O}_\prec(\Psi^2)\,,
\,,
}
in terms of the small control parameter $ \Psi \leq N^{-c}$. 
Writing $ \brm{w} = \brm{T}\1(\brm{v}-(\brm{e},\brm{v})\2\brm{e}) $, and recalling $ [\brm{T},\brm{S}] = \brm{0} $ (actually we need only the case $ \brm{T} =  \brm{S} $ here) and $ (\1\brm{f},\brm{v}) = 0 $ by \eqref{f-projection of v = 0}, we obtain:
\bels{f-projection of w = 0}{
(\1\brm{f},\brm{w}) \,=\, 0
\,. 
}   
Thus by expressing \eqref{SCE for FA} in the vector form, 
\bels{linearised SCE for w}{
(1-m^2\brm{S}\1)\1\brm{w} \,=\, \mcl{O}_{\msp{-1}\prec}\msp{-1}(\Psi^2)
\,,
}
we see that $ 1-m^2\brm{S} $ can be also inverted in the subspace orthogonal to both the $ +1 $ and $ -1 $ eigendirections. Hence \eqref{linearised SCE for w} yields
\[
\brm{w} \,=\, \mcl{O}_{\msp{-2}\prec}\msp{-1}\bigl(\2\widehat{\Gamma}\,\Psi^2\bigr)
\,,
\] 
which is exactly the fluctuation averaging bound (4.14) of \cite{Erdos13-LoclScGenRMT} with $ \wti{\Gamma} $ updated to $ \wht{\Gamma} $. 

Besides these observations and the replacement of $ \wti{\Gamma} $ by $ \wht{\Gamma} $ the proof from \cite{Erdos13-LoclScGenRMT} can be carried out without further modifications. 
 
%%% References %%%%%%%%%%%%%%%%%%%%%%%%%%%%%%%%%%%%%

\end{document}